\documentclass[11pt]{article}
\usepackage{amssymb,amsmath,bm}
\usepackage[colorlinks=true, urlcolor=blue,linkcolor=blue, citecolor=blue]{hyperref}
\usepackage{mathrsfs,verbatim}
\usepackage{caption,float}
\usepackage{constants,color}
\usepackage{enumerate}
\usepackage{bbm}
\usepackage{appendix}
\usepackage{graphicx, caption,subcaption}
\begin{document}
\newcommand{\bea}{\begin{eqnarray}}
\newcommand{\ena}{\end{eqnarray}}
\newcommand{\beas}{\begin{eqnarray*}}
\newcommand{\enas}{\end{eqnarray*}}
\newcommand{\beq}{\begin{equation}}
\newcommand{\enq}{\end{equation}}
\def\qed{\hfill \mbox{\rule{0.5em}{0.5em}}}
\newcommand{\bbox}{\hfill $\Box$}
\newcommand{\ignore}[1]{}
\newcommand{\ignorex}[1]{#1}
\newcommand{\wtilde}[1]{\widetilde{#1}}
\newcommand{\mq}[1]{\mbox{#1}\quad}
\newcommand{\bs}[1]{\boldsymbol{#1}}
\newcommand{\qmq}[1]{\quad\mbox{#1}\quad}
\newcommand{\qm}[1]{\quad\mbox{#1}}
\newcommand{\nn}{\nonumber}
\newcommand{\Bvert}{\left\vert\vphantom{\frac{1}{1}}\right.}
\newcommand{\To}{\rightarrow}
\newcommand{\supp}{\mbox{supp}}
\newcommand{\law}{{\cal L}}
\newcommand{\D}{{\cal D}}
\newcommand{\Z}{\mathbb{Z}}
\newcommand{\E}{\mathbb{E}}
\newcommand{\Pro}{\mathbb{P}}
\newcommand{\eff}{\mathrm{Eff}}
\newcommand{\dist}{\mathrm{dist}}


\newcommand{\var}{\mathrm{Var}}
\newcommand{\cov}{\mathrm{Cov}}

\newcommand{\ucolor}[1]{\textcolor{blue}{#1}}  
\newcommand{\ucomm}[1]{\marginpar{\tiny\ucolor{#1}}}  

\newcommand{\bcolor}[1]{\textcolor{red}{#1}}  
\newcommand{\bcomm}[1]{\marginpar{\tiny\icolor{#1}}}  

\newtheorem{theorem}{Theorem}[section]
\newtheorem{corollary}{Corollary}[section]
\newtheorem{conjecture}{Conjecture}[section]
\newtheorem{proposition}{Proposition}[section]
\newtheorem{lemma}{Lemma}[section]
\newtheorem{definition}{Definition}[section]
\newtheorem{example}{Example}[section]
\newtheorem{remark}{Remark}[section]
\newtheorem{case}{Case}[section]
\newtheorem{condition}{Condition}[section]
\newcommand{\pf}{\noindent {\it Proof:} }
\newcommand{\proof}{\noindent {\it Proof:} }
\frenchspacing


%

\title{\bf  Clique structure and other network properties of the tensor product of  Erd\H{o}s-R\'enyi graphs} 
\author{\"{U}m\.{i}t I\c{s}lak\footnote{Bo\u{g}azi\c{c}i University, Mathematics Department, Istanbul, Turkey. email: umit.islak1@bogazici.edu.tr}\hspace{0.1in}\footnote{Middle East Technical University, Institute of Applied Mathematics, Ankara, Turkey. email: uislak@metu.edu.tr} \hspace{0.2in}  \quad Buğra İncekara\footnote{Technical University of Munich, Department of Mathematics, Munich, Germany. email: bugra.incekara@tum.de
}}\vspace{0.25in}
\maketitle
\begin{abstract}
We analyze the number of cliques of given size and the size of the largest clique in tensor product $G \times H$ of two Erd\H{o}s-R\'enyi graphs $G$ and $H$. Then an extended clustering coefficient is introduced and is studied for $G \times H$. Restriction to the standard clustering coefficient has a direct relation to the local efficiency of the graph, and the results are also interpreted in terms of the efficiency. As a last statistic of interest, the number of isolated vertices is  analyzed for $G \times H$. The paper is concluded with a discussion of the modular product of random graphs, and the relation to the maximum common subgraph problem.
\end{abstract}

\section{Introduction}

An Erd\H{o}s-R\'enyi graph is  a random graph in which each possible edge is present with some fixed probability $p \in (0,1)$, and the presences of edges are independent of each other. Since the pioneering work (\cite{er:59}, \cite{er:60}) of Erd\H{o}s and R\'enyi in late 50's, there had been a tremendous development in the field of random graphs, and much more realistic models had been introduced and studied in the literature. See \cite{bol:98} and \cite{hofstadt:24} for general overviews. Real life networks had been analyzed in terms of their properties such as the clustering and efficiency, and these have had important applied and theoretical reflections. In this note, we go back to the setting of  Erd\H{o}s-R\'enyi graphs and analyze the clique structure and network properties of the tensor product of such two random graphs. Our interest in this problem  stems from the relation between the closely related  modular product of graphs and the maximum common subgraph problem, and relevant discussions will be included below.

Graph products in general appear in several distinct contexts and applications, and various graph products are well-studied in the literature. The reader is referred to  \cite{west:01} for graph theoretical foundations, and \cite{hammack:11}  for a detailed background on graph products.  As already noted, our interest here will be on the tensor product and modular product of graphs. Focusing on the former, given two graphs $G= (V_1, E_1)$ and $H = (V_2, E_2)$, the tensor product $G \times H$ is the graph with vertex set $\{(v_1, v_2): v_1 \in V_1, v_2 \in V_2  \}$  in which two vertices  $(v_1, v_2)$ and $(v_1', v_2')$  are adjacent if and only  if $v_1$ is adjacent to $v_1'$ in $G$ and  $v_2$ is adjacent to $v_2'$ in $H$. On the other hand, the modular product of $G \times_m H$ has the same vertex   set with $G \times H$, but two vertices $(v_1, v_2)$ and $(v_1', v_2')$  are adjacent if and only  if  one of the following two is true: (i) $v_1$ is adjacent to $v_1'$ in $G$ and  $v_2$ is adjacent to $v_2'$ in $H$, or, (ii) $v_1$ is not adjacent to $v_1'$ in $G$ and  $v_2$ is not adjacent to $v_2'$ in $H$.

A motivation for the tensor product of graphs can be given as follows. Suppose that the vertices of the first and second graphs correspond to certain sets of men and women, respectively. Using the friendships among men, we may form a graph $G_1$, and similarly via the friendships in $G_2$ another graph $G_2$. Then the tensor product  $G_1 \times G_2$ consists of the pairs of a man and a woman, and the corresponding relation in product graph is based on whether men and women are friends within their group. A second exemplary case could consist of two sets of missions. The first graph $G_1$ has the first set of missions as its vertex set, and two missions are adjacent if they could be run simultaneously. The second graph is  formed similarly. The resulting tensor product then holds the information of pairs of missions from the two sets as vertices, and two of the vertices are adjacent if the two missions corresponding to both first and second one can be run simultaneously.

 In this manuscript,
we are mainly interested in the clique structure of the tensor product of two Erd\H{o}s-R\'enyi graphs, but we also look at certain other statistics   in the complex networks literature. In particular, besides the number of cliques of a given size, we will be interested in the related problems of analyzing  the clustering coefficient and efficiency. The latter topics in networks had recently been studied for distinct graph products; see, for example, \cite{agj:23} and \cite{dme:22}. Yet another statistic we will analyze is the number of isolated vertices in tensor product of random graphs, and this is not directly related to the underlying clique structure of the graphs. 

Before continuing with the content of the present manuscript, let us note that the main motivation of our study arises actually from another graph product, namely, the modular product. This special  product has a relation to comparison and similarity of two graphs, which is an important topic in  both theory and applications. See \cite{wm:20} for a comprehensive survey on graph comparisons. One well known notion used for such purposes is the maximum common subgraph which  is a graph that is an induced subgraph of the compared graphs, and that has the maximum possible vertices. It is well known that the size of maximum common subgraph  of two graphs is the same as the clique number of the modular product  of $G$ and $H$ (\cite{bb:76}). Our initial purpose in this work was to study the similarity of two independent Erd\H{o}s-R\'enyi graphs through an analysis of the modular  product of these graphs. However, we were not currently able to do a through analysis of the modular product case, and we present our results on the tensor products in this manuscript.  That said, some results obtained here are extendable to the modular product setup, and we plan to do a detailed analysis in an upcoming work. 
 
Let us now fix some notation.  $F, G,$ and $H$ are used for graphs. The graphs considered are all assumed to be simple. $V$ and $E$ denote the vertex sets and edge sets, respectively. For given graphs $G,H$, $G \times H$ and $G \times_m H$ correspond to the tensor and modular products of $G$ and $H$, respectively.   $d^G(v)$ is the degree of vertex $v \in V$ in graph $G$. For  $v \in V$,  $N(v)$ is the set of vertices that are adjacent to $v$. In particular $v \notin N(v)$. $I_A$ is used for the indicator function of the event $A$. We write $\log$ instead of $\log_2$. Other notation will be introduced below in its first appearance. 

The rest of the paper is organized as follows. Let $G$ and $H$ be two independent Erd\H{o}s-R\'enyi graphs with parameters $n \in \mathbb{N}$ and $p \in (0,1)$, and $G \times H$ be their tensor product. Next section obtains results on  the expected number of $k$-cliques in $G \times H$, which are to be later used for understanding the largest clique and the extended clustering coefficient. Section \ref{sec:maxclique} shows that the largest clique in $G \times H$ has size approximately $2 \log_{1/p} n$ as in the case of a single Erd\H{o}s-R\'enyi graph. Section \ref{sec:extcc} introduces an extended notion of clustering coefficient, and proves a strong law of large numbers for the corresponding coefficient in $G \times H$. Section \ref{sec:isol} studies the expected number of isolated vertices in $G \times H$, and complements this with a weak law of large numbers in a special case. The paper is concluded in Section \ref{sec:conc} with a discussion of possible questions for future work. In particular, the largest clique in modular product of Erd\H{o}s-R\'enyi graphs and its relation to the maximum common subgraph problem is questioned. 

\section{The number of $k$-cliques}\label{sec:numofkcliques}

The purpose of this section is to study the number of $k$-cliques in tensor product of  Erd\H{o}s-R\'enyi graphs, which in turn later will be used for analyzing the size of the corresponding largest clique. The latter problem in a single Erd\H{o}s-R\'enyi  graph  is first studied in \cite{matula:76}. Our related proofs below will be following \cite{as:16} who sketch an argument for the symmetric case. We will extend their proof to the non-symmetric case with detailed calculations which will allow us to obtain explicit variance upper bounds for the number of $k$-cliques when $k$ is  smaller than $2 \log_{1/p} n - M \log_{1/p} \log n$ where $M > 4$. We will also restrict our observations to the fixed $k$ case, and will use the corresponding work in our analysis of the  extended clustering coefficient below.  

Now, towards understanding the cliques in tensor product of graphs, first observe that a given $k$-clique $\{u_1,\ldots,u_k\}$ of $G$, and a  $k$-clique $\{v_1,\ldots,v_k\}$ of $H$ provide $k!$-cliques in the tensor product.  $\{(u_1,v_1),(u_2,v_2),\ldots,(u_k,v_k)\}$ is one of them, and keeping the role of $u$'s the same, but reordering $v$'s we obtain $k!$ corresponding cliques.)  
 
 On the other hand, given a $k$-clique $\{(u_1,v_1),(u_2,v_2),\ldots,(u_k,v_k)\}$ of the tensor product $G \times H$, $\{u_1,\ldots,u_k\}$ and   $\{v_1,\ldots,v_k\}$ clearly forms a clique of $G$ and $H$, respectively. Therefore we have the following result.
 
\begin{proposition}
    Let $G$ and $H$ be two graphs on $n$ vertices. Denote the number of $k$-cliques in $G$ and $H$ by $X_k$ and $Y_k$, respectively. Then  the number of $k$-cliques $Z_k$ in the tensor product $G \times H$ is given by $$Z_k = k! X_k Y_k.$$
 \end{proposition}

Now, given the proposition we have the following useful observation  for random graphs case. 

\begin{proposition}\label{propn:firsts}
    Let $G$ and $H$ be two independent  Erd\H{o}s-R\'enyi graphs on $n$ vertices with attachment probability $p$. Let $m \in \mathbb{N}$. Denoting the number of $k$-cliques in $G$, $H$ and $G \times H$ by $X_k$, $Y_k$ and $Z_k$, respectively, we have the following relation: 
    $$\mathbb{E}[Z_k^m] = (k!)^m \left( \mathbb{E}[X_k^m] \right)^2. $$ 
\end{proposition}

\textbf{Proof.}  We have 
$$\mathbb{E}[Z_k^m] = \mathbb{E}[(k! X_k Y_k))^m] = (k!)^m \mathbb{E}[X_k^m] \mathbb{E}[Y_k^m] = (k!)^m \left( \mathbb{E}[X_k^m] \right)^2$$
 \hfill $\square$

Specializing on $m=1,2$ cases give the following. 

\begin{theorem}
    Let $G$ and $H$ be two independent  Erd\H{o}s-R\'enyi graphs on $n$ vertices with attachment probability $p$.  Let $Z_k$ be  the number of $k$-cliques in $G \times H$. Then, we have 
\begin{equation}\label{eqn:meanzk}
    \mathbb{E}[Z_k] =  k! \binom{n}{k}^2 p^{k(k-1)},
\end{equation}
     and 
     \begin{equation}\label{eqn:varzk}
     \var(Z_k) = (k!)^2 ((\var(X_k))^2 + 2 \var(X_k) (\mathbb{E}[X_k])^2),
     \end{equation} where $X_k$ is the number of $k$-cliques in $G$.
\end{theorem}
\textbf{Proof.}  Using Proposition \ref{propn:firsts}, we have 
$$
\mathbb{E}[Z_k] = k! (\mathbb{E}[X_k])^2 = k! \left( \binom{n}{k} p^{\binom{k}{2}}\right)^2,$$ from which the first observation follows. For the second,  again using Proposition \ref{propn:firsts},
\begin{eqnarray*}
\var(Z_k) &=& \mathbb{E}[Z_k^2] - (\mathbb{E}[Z_k] )^2 = (k!)^2 (\mathbb{E}[X_k^2])^2 - ((k!) (\mathbb{E}[X_k])^2)^2 = (k!)^2 ((\mathbb{E}[X_k^2])^2  - (\mathbb{E}[X_k])^4)    \\
&=&  (k!)^2 ((\var(X_k) + (\mathbb{E}[X_k])^2)^2 - (\mathbb{E} [X_k])^4)) = (k!)^2  ((\var(X_k))^2 + 2 \var(X_k) (\mathbb{E}[X_k])^2).
\end{eqnarray*}
 \hfill $\square$

 The second part expresses the variance of the number of cliques in the tensor product in terms of the corresponding variance in an individual Erd\H{o}s-R\'enyi graph. The latter is previously studied in Theorem 4.5.1 of \cite{as:16} for $p=1/2$ in a close framework. The following  upper bound on $\var(X_k)$ is obtained similarly to the steps given in  \cite{as:16} for $p =1/2$. We include the proof for the general $p$ case in the appendix. Let us again note that  the  argument given in the proof will also be useful later in our discussion on the generalized clustering coefficient. 

 \begin{proposition}\label{propn:basicupbound}
     Let $X_k$ be the number of $k$-cliques in an Erd\H{o}s-R\'enyi graph $G$ on $n$ vertices with attachment probability $p$. Then, 
     $$\var(X_k) \leq  \binom{n}{k}  p^{\binom{k}{2}} + \binom{n}{k}  p^{\binom{k}{2}} \sum_{i=2}^{k-1} \binom{k}{i} \binom{n  - k}{k - i} p^{\binom{k}{2} - \binom{i}{2}}.$$
 \end{proposition}

The proposition has the following corollary whose proof is again left to the appendix.

\begin{corollary}\label{cor:varbound}
(i) Let $M > 4$. If $k =  2 \log_{1/p} n - M \log_{1/p} \log n$, then for sufficiently large $n$,  we have 
\begin{equation}\label{eqn:varmainbd}
\var(X_k) \leq \binom{n}{k}  p^{\binom{k}{2}} + \binom{n}{k} \binom{n-k}{k-2}  p^{k(k-1)} \frac{k^3}{2} p^{-1}.    
\end{equation}

(ii) The  bound in \eqref{eqn:varmainbd} holds also for the case where $k$ is fixed. 
 \end{corollary}

Indeed as the proof for the corollary will show, the same variance bound will hold when $k$ is smaller than $k =  2 \log_{1/p} n - M \log_{1/p} \log n$ where $M > 4$. For fixed $k$ case, we also note the following.

\begin{corollary}\label{cor:2}
    When $k$ is fixed, we have, $$\var(X_k) \leq C n^{2k - 2},$$ and $$\var(Z_k) \leq D n^{4k -2},$$ where $C, D$ are constants depending on $p$ and $k$, but not on $n$. 
\end{corollary}

\textbf{Proof.} 
The first statement is clear from \eqref{eqn:varmainbd}. For the second one we recall      $$\var(Z_k) = (k!)^2 ((\var(X_k))^2 + 2 \var(X_k) (\mathbb{E}[X_k])^2)$$ from  \eqref{eqn:varzk}, and  use the estimates $\var(X_k) \leq C n^{2k - 2},$ and  $\mathbb{E}[X_k] \leq n^k$. 
\hfill $\square$ 

Before concluding this section, let us note that the variance bounds we have found can be  inserted in \eqref{eqn:varzk} for growing $k$ (such as $k =  2 \log_{1/p} n - M \log_{1/p} \log n$)   in order to obtain  similar bounds for the variance $\var(Z_k)$ corresponding to the tensor product. The resulting expressions turn out to be  complicated, and will not be stated since they will not be useful in that form below.

\section{The clique number of the tensor product}\label{sec:maxclique}

The purpose of this section is to prove two results on the clique number (i.e., the number of vertices in a maximum clique) of the tensor product of two independent Erd\H{o}s-R\'enyi graphs.   The conclusion for $p=1/2$ case will be that the clique number behaves like $2 \log n$ as the number of vertices in both vertices tend to infinity. So in this case  the clique number of the tensor product behaves the same as the clique number of an individual Erd\H{o}s-R\'enyi graph. 

\begin{theorem}
     Let $G$ and $H$ be two independent  Erd\H{o}s-R\'enyi graphs on $n$ vertices with attachment probability $p$. Let $Z_k$ be the number of $k$-cliques in $G \times H.$ Then, there exists some non-negative sequence $a_n$ with $a_n = o(1)$, such that if  $k^* = (2 + a_n) \log_{1/p} n$, then $$\Pro(Z_{k^*} > 0) \rightarrow 0, \quad \text{as } n \rightarrow \infty.$$
\end{theorem}

\textbf{Proof.} Recalling the mean of $Z_{k}$ from \eqref{eqn:meanzk},
\begin{eqnarray*}
    \mathbb{E}[Z_k] &=& k! \binom{n}{k}^2 p^{k(k-1)} = k! \frac{(n (n-1)\cdots (n - k +1))^2}{(k!)^2} p^{k(k-1)} \leq \frac{n^{2k}}{k! }p^{k(k-1)} \\
    &\leq& n^{2k } p^{k(k-1)} = (1/p)^{2k \log_{1/p} n - k (k-1)} = (1/p)^{k (2 \log_{1/p}n - k + 1)}.
\end{eqnarray*}
Now choosing $a_n = \frac{\log_{1/p} \log n}{\log_{1/p}  n}$, and using $k^* = (2 + a_n) \log_{1/p}  n$, we obtain 
$$\mathbb{E}[Z_{k^*}] \rightarrow 0, \quad \text{as } n \rightarrow \infty.$$ 
Thus, using Markov's inequality,
$$\mathbb{P}(Z_{k^*} > 0) = \mathbb{P}(Z_{k^*} \geq 1) \leq \E[Z_{k^*}] \rightarrow 0, \quad \text{as } n \rightarrow \infty,$$ as asserted.  \hfill $\square$

\begin{theorem}
     Let $G$ and $H$ be two independent  Erd\H{o}s-R\'enyi graphs on $n$ vertices with attachment probability $p$. Let $Z_k$ be the number of $k$-cliques in $G \times H.$ Then, there exists some non-negative sequence $a_n$ with $a_n = o(1)$, such that if  $k_* = (2 - a_n) \log_{1/p} n$, then $$\Pro(Z_{k_*} > 0) \rightarrow 1, \quad \text{as } n \rightarrow \infty.$$
\end{theorem}
\textbf{Proof.}  Let $M > 4$ and $k_* = (2 - a_n) \log_{1/p} n$ with $a_n = \frac{M \log_{1/p} \log n}{\log_{1/p} n}$. Note that proving $\Pro(Z_{k^*} > 0) \rightarrow 1$ is equivalent to $\Pro(Z_{k^*} = 0) \rightarrow 0$. Recall the inequality, $$\Pro(Z_{k^*} = 0)\leq \frac{\var(Z_{k_*})}{(\E[Z_{k_*}])^2}.$$
Using \eqref{eqn:meanzk} and \eqref{eqn:varzk}, we have 
$$\Pro(Z_{k_*} = 0)\leq \frac{\var(Z_{k_*})}{(\E[Z_{k_*}])^2} = \frac{(k_*!)^2 ((\var(X_{k_*}))^2 + 2 \var(X_{k_*}) (\mathbb{E}[X_{k_*}])^2)}{(k_*!)^2 (\E[X_{k_*}])^4} = \left(\frac{\var(X_{k_*})}{(\E[X_{k_*}])^2} \right)^2  + \frac{2 \var (X_{k_*})}{(\E[X_{k_*}])^2}.$$ Hence we are just left with showing that $\frac{\var(X_{k_*})}{(\E[X_{k_*}])^2} \rightarrow 0$ as $n \rightarrow \infty$.
Using Corollary \ref{cor:varbound}, it will be enough to show that 
$$\frac{ \binom{n}{k_*}  p^{\binom{k_*}{2}} + \binom{n}{k_*} \binom{n-k_*}{k_*-2}  p^{k_*(k_*-1)} \frac{k_*^3}{2} p^{-1}}{\binom{n}{k_*}^2 p^{k_* (k_*-1)}} \rightarrow 0, \quad \text{as } n \rightarrow \infty.$$

For this purpose we will show, $\frac{1}{\binom{n}{k_*}p^{\binom{k_*}{2}}} \rightarrow 0$, and $\frac{\binom{n - k_*}{k_*-2} k_*^3}{\binom{n}{k_*}} \rightarrow 0$. Beginning with the latter, observe that 
$$\frac{\binom{n - k_*}{k_*-2} k_*^3}{\binom{n}{k_*}} = \left(\frac{n-k_*}{n} \frac{n-k_*-1}{n-1} \cdots \frac{n-2k_*+3}{n-k_*+3}\right)\left( \frac{1}{n -k_*+2} \frac{1}{n -k_*+1} \right) \frac{k_*!}{(k_*-2)!}  k_*^3 \leq \frac{k_*^5}{(n - k_* +2)^2}.$$ With the given selection of $k_*$, this right-hand side clearly converges to 0. Next, in order to show  $\frac{1}{\binom{n}{k_*}p^{\binom{k_*}{2}}} \rightarrow 0$, it is enough to see $\binom{n}{k_*}p^{\binom{k_*}{2}} \rightarrow \infty$. 
For this purpose, writing $r = 1/p$, and substituting $k_* = 2  \log_{1/p} n - M \log_{1/p} \log n$, we observe that 
\begin{eqnarray*}
\binom{n}{k_*}p^{\binom{k_*}{2}} &\geq& \left( \frac{n}{k_*} \right)^{k_*} p^{\frac{k_*(k_*-1)}{2}} = p^{k_* \log_p n - k_* \log_p k_* + \frac{k_*(k_*-1)}{2} }  = p^{\frac{k_*}{2} (2 \log_p n - 2 \log_p k_* + k_* - 1)} \\
&=& p^{\frac{k_*}{2} (-M\log_r \log n +2 \log_r (2 \log_r n - M\log_r \log n) - 1)} \geq  p^{\frac{k_*}{2} (-M\log_r \log n +2 \log_r (\log_r n) - 1)}. 
\end{eqnarray*}
We see that the right most term  tends to infinity, and thus the proof is done.  \hfill $\square$

 \section{An extended clustering coefficient}\label{sec:extcc} 

\subsection{Definition and basic observations}

Let $F$ be a graph.  The \textit{standard  clustering coefficient} of a  given vertex $u \in V(F)$ is defined to be $$CC^F(u) = \frac{\sum_{v,w \in N(u)} \mathbf{1}(v \leftrightarrow w)}{\binom{d^F(u)}{2}}.$$ Although the underlying notion was studied earlier, the concept of  clustering coefficient got   particular attention after the seminal paper of Watts and Strogatz on the small world model \cite{ws:98}. Since then there had been several variations and generalizations, some of which can be found in \cite{fronczaketal:02} and \cite{haoetal:18}. 

Focusing on the standard clustering coefficient $CC^F(u)$, the clustering coefficient of vertices in tensor products of graphs were recently studied in \cite{dme:22}. In this paper, the main result is based on writing the number of triangles incident to a vertex $(u_1,v_1)$ in tensor product $G \times H$ in terms of the number of  triangles in $G$ and $H$. 

To give a quick look at their main observation, clearly, every triangle in $G$ containing $u_1$, and  every triangle in $H$ containing $v_1$ will yield a triangle in the product (and also the other way). The question is how many of them? Since one corner is already $(u_1, v_1)$, there are two options for the others   $(u_2, v_2)$, $(u_3, v_3)$ or $(u_2, v_3)$, $(u_3, v_2)$. Basically we consider $(3 - 1)! = 2!$ permutations for the remaining ones, and thus there are two of them. This observation is easily extendable to a more general framework which we do next. 

  Let $N(w)$ denote the neighborhood of a vertex $w$ (again $w$ not included). Then we define the $k$-\textit{clustering coefficient} of $W$ to be $$C_k^G(w) = \frac{\#\text{of } k-\text{cliques in } N(w)}{\binom{|N(w)|}{k}} =: \frac{A_k^G(w)}{\binom{|N(w)|}{k}}.$$ Thus, $A_k^G(w)$ is the number of $k$-cliques incident to the vertex $w$.  When $k = 2$, this definition reduces to  the standard clustering coefficient. The definition we give is related, but not exactly the same as the one previously introduced in \cite{haoetal:18}. 

Now, the question is whether a similar formulation to \cite{dme:22} is possible for $A_k$ in tensor products. 

\begin{proposition}\label{propn:tridecomp}
Let $G, H$ be graphs. 
Let $u_1 \in V(G)$ and $v_1 \in V(H)$. Then 
$$A_k^{G \times H}((u_1,v_1)) = k! A_k^G(u_1) A_k^H(v_1).$$   
\end{proposition}

\textbf{Proof.}
For $u_1 \in V(G)$ and $v_1 \in V(H)$, we are interested in the number of $k$-cliques  that are incident to $(u_1,v_1) $ in $G \times H$. As we discussed before for $k=2$ case, every $(k+1)$-clique in $G$ containing $u_1$, and  every  $(k+1)$-clique in $H$ containing $v_1$ will yield a  $(k+1)$-clique in the product (and also the other way). How many $(k+1)$-cliques are there in $G \times H$ for the given pair of $(k+1)$-cliques in $G$ and $H$? Since one corner is already $(u_1, v_1)$, there are   $k!$ options for the other ones.   
Hence we conclude, 
$$A_k^{G \times H}((u_1,v_1)) = k! A_k^G(u_1) A_k^H(v_1).$$   
\hfill $\square$

Next  result will be useful when we analyze the clustering coefficient of the tensor product in the next subsection.

\begin{proposition}
    Let $u \in V(G)$ and $v \in V(H)$. Then,
    $$C_k^{G \times H}((u,v)) =   k!  C_k^G(u) C_k^H(v) D(G, H, k, u, v),$$ where $$D(G, H, k, u, v) = \frac{\binom{d^G(u)}{k} \binom{d^H(v)}{k}}{ \binom{d^G(u)d^H(v)}{k}}.$$
\end{proposition}

\textbf{Proof.} We recall  that for an individual graph $F$, $$C_k^F(w) =     \frac{A_k^F(w)}{\binom{d^F(w)}{2}}.$$ 
Now let $G$  and $H$ be two graphs on $n$ vertices as given in the theorem. Let $u \in V(G)$ and $v \in V(H)$.   Observe that we have $$C_k^{G \times H}((u,v)) = \frac{A_k^{G \times H}((u,v))}{\binom{d^{G \times H}((u,v))}{k}} = \frac{k ! A_k^G(u) A_k^H(v)}{ \binom{d^G(u)d^H(v)}{k}}.$$  Here the last step uses our observation from the Proposition \ref{propn:tridecomp}, along with the  fact that $$d^{G \times H}((u,v)) = d^G(u)d^H(v),$$ which can be verified  by using the definition of tensor product. Rewriting this last expression we obtain 
 $$C_k^{G \times H}((u,v))  = \frac{k ! A_k^G(u) A_k^H(v)}{ \binom{d^G(u)d^H(v)}{k}} = k! \frac{A_k^G(u) A_k^H(v)}{\binom{d^G(u)}{k} \binom{d^H(v)}{k}} \frac{\binom{d^G(u)}{k} \binom{d^H(v)}{k}}{ \binom{d^G(u)d^H(v)}{k}} = k!  C_k^G(u) C_k^H(v) D(G, H, k, u, v),$$  where we set 
 $D(G, H, k, u, v) =\frac{\binom{d^G(u)}{k} \binom{d^H(v)}{k}}{ \binom{d^G(u)d^H(v)}{k}},$  as in the statement of the proposition.
\hfill $\square$

 \subsection{Asymptotic analysis of $C_k^{G\times H}((u,v))$ for Erd\H{o}s-R\'enyi graphs}
 
Our main result will be  a strong law of large numbers  $C_k^{G\times H}((u,v))$ when $G$ and $H$ are independent Erd\H{o}s-R\'enyi graphs with parameters $n,p$, and when $n \rightarrow \infty$. This first requires the following analogous result for a single random graph. 

\begin{proposition}\label{propn:asconv}
Let $G$  be  an Erd\H{o}s-R\'enyi  graph with parameters $n\in \mathbb{N}$ and $p \in (0,1)$. Let $k \in \mathbb{N}$ be fixed. Then, 
$$C_k^G (u) \rightarrow p^{\binom{k}{2}},  \quad n \rightarrow  \infty,$$ where the convergence is with probability one. 
\end{proposition}

\textbf{Proof.} 
Writing $$C_k^G(u) = \frac{A_k^G(u)}{\binom{d^G(u)}{k}} = \frac{A_k^G(u) / n^k}{\binom{d^G(u)}{k} /n^k}, $$ it suffices to show that  
\begin{itemize}
    \item[(i)] $\binom{d^G(u)}{k} /n^k \rightarrow \frac{p^k}{k!}$, and,
    \item[(ii)] $A_k^G(u) / n^k \rightarrow \frac{p^{\binom{k + 1}{2}}}{k!}$,
\end{itemize}
 where the convergences are with probability one, as $n \rightarrow \infty$. For the first assertion, since $k$ is fixed, observe that $$ 
\frac{\binom{d^G(u)}{k}}{n^k}  = \frac{1}{k!} \left( \frac{d^G(u)}{n} \right) \left( \frac{d^G(u) - 1}{n} \right) \cdots \left( \frac{d^G(u) - k + 1}{n} \right) \rightarrow \frac{1}{k!} p \cdot p \cdots p = \frac{p^k}{k!},$$ with probability one. Here, we use the fact that $d^G(u)$ is a sum of independent and identically distributed Bernoulli random variables, and  the standard  strong law of large numbers. 

Regarding the second claim we observe that 
$$\frac{A_k^G(u)}{n^k} = \frac{\sum_{S \subset \{1,2,\ldots,n\} \setminus \{u\},  \#S = k}  \mathbf{1}(S \cup \{u\} \text{ forms a  clique})}{n^k}.$$ Denoting the summation as $\sum_S$, the numerator  can be rewritten as $\sum_S \mathbf{1}(S \text{ forms a clique}) \mathbf{1}(u \leftrightarrow S)$. 
Now, writing $I_S = \mathbf{1}(S \text{ forms a clique})$, we then have
\begin{eqnarray*}
\var(A_k^G(u)) &=& \var \left(\sum_S I_S \mathbf{1} (u \leftrightarrow S) \right)     \\
&\leq& \sum_S \mathbb{E}[I_S \mathbf{1} (u \leftrightarrow S)] + \sum_{S \neq T} \mathbb{E}[I_S I_T \mathbf{1} (u \leftrightarrow S, u \leftrightarrow T)] \\
&\leq&  \sum_S \mathbb{E}[I_S] + \sum_{S \neq T} \mathbb{E}[I_S I_T]. 
\end{eqnarray*}

But then following the same steps in the proof of Proposition  \ref{propn:basicupbound}   with replacing $n$ by $n - 1$, we arrive at the conclusion $$\var(A_k^G(u)) \leq L n^{2k-2}$$ for some constant $L$ as in   Corollary \ref{cor:2}. Then, using Chebyshev's inequality, we have $$\mathbb{P} \left( \left|\frac{A_k^G(u)}{n^k} -\mathbb{E} \left[ \frac{A_k^G(u)}{n^k} \right] \right| > \epsilon \right) \leq \frac{\var(A_k^G(u))}{\epsilon^2 n^{2k}} \leq \frac{L n^{2k-2}}{\epsilon^2 n^{2k}} = \frac{L}{\epsilon^2 n^2}.$$
Adding through all possible $n$'s and using Borel-Cantelli's first lemma we conclude that 
$\frac{A_k^G(u)}{n^k} \rightarrow \mathbb{E} \left[ \frac{A_k^G(u)}{n^k} \right]$ with probability one. Since $\mathbb{E} \left[ \frac{A_k^G(u)}{n^k} \right] \sim \frac{p^{\binom{k+1}{2}}}{k!}$, our claim that $A_k^G(u) / n^k \rightarrow \frac{p^{\binom{k + 1}{2}}}{k!}$ with probability one follows. 

Combining these observations, we conclude $$C_k^G(u)  = \frac{A_k^G(u) / n^k}{\binom{d^G(u)}{k} /n^k} \rightarrow \frac{\frac{p^{\binom{k + 1}{2}}}{k!}}{\frac{p^k}{k!}} = p^{\binom{k}{2}},$$  with probability  as required.
 \hfill $\square$
 
Now we are ready to state the main result of this section  on extended clustering of  tensor product of  random graphs. 

\begin{theorem}\label{thm:asconvofck}
Let $G$  and $H$ be two independent Erd\H{o}s-R\'enyi  graphs, 
each with parameters $n\in \mathbb{N}$ and $p \in (0,1)$. Let $u \in V(G) $ and $v \in V(H)$.  We have $C_k^{G \times H} ((u,v))\longrightarrow  p^{2 \binom{k}{2}} $ with probability 1 as $n \rightarrow \infty.$
\end{theorem}

\textbf{Proof.} 
Recall from Proposition that \ref{propn:tridecomp}  $$C_k^{G \times H}((u,v)) =   k!  C_k^G(u) C_k^H(v) D(G, H, k, u, v),$$  where
$$D(G, H, k, u, v) =\frac{\binom{d^G(u)}{k} \binom{d^H(v)}{k}}{ \binom{d^G(u)d^H(v)}{k}} .$$   
We already know by Proposition \ref{propn:asconv} that  both $C_k^G(u)$ and $C_k^H(v)$ converge to $p^{\binom{k}{2}}$ with probability one. Hence   it will suffice to show 
   $$D(G, H, k, u, v) \rightarrow \frac{1}{k!}, $$  with probability 1 as $n\rightarrow \infty$.
 
 In order to prove this, let us first write 
$$D(G, H, k, u, v)  = \frac{\binom{d^G(u)}{k} \binom{d^H(v)}{k}}{ \binom{d^G(u)d^H(v)}{k}} = \frac{1}{k!} \frac{(d^G(u))_k (d^H(v))_k}{(d^G(u) d^H (v))_k}.$$ (Here $(a)_k$  is the $k$th falling factorial of $a$.) So we would like to conclude that $\frac{(d^G(u))_k (d^H(v))_k}{(d^G(u) d^H (v))_k} \rightarrow 1$. For this purpose we write 

\begin{eqnarray*}
\frac{(d^G(u))_k (d^H(v))_k}{(d^G(u) d^H (v))_k} &=& \prod_{i = 0}^{k-1}  \frac{(d^G(u) - i) (d^H(v) - i)}{(d^G(u)d^H(v) - i)}  \\ &=&  \prod_{i = 0}^{k-1}  \left( \frac{d^G(u) d^H(v) - i -i(d^G(u) + d^H(v)) + i^2 + i}{d^G(u) d^H(v) - i }  \right). 
\end{eqnarray*}
Now for each fixed $i$, looking at, $$\frac{d^G(u) d^H(v) - i -i(d^G(u) + d^H(v)) + i^2 + i}{d^G(u) d^H(v) - i }  = 1 - \frac{i (d^G(u) +  d^H(v))}{d^G(u) d^H(v) - i} + \frac{i^2 + i}{d^G(u) d^H(v) - i},$$ we see that the last two terms tend to zero with probability one  as $n \rightarrow \infty$, since both $d^G(u)$ and $ d^H(v)$ tend to infinity as $n \rightarrow \infty$, again with probability one. 
Combining these observations, we conclude that $D(G, H, k, u, v) \rightarrow \frac{1}{k!}$ as required. 

Hence, the proof that  $C_k^{G \times H} ((u,v))\longrightarrow  p^{2 \binom{k}{2}} $ with probability 1, as stated in the theorem is done. \hfill $\square$

We close this section by noting that the theorem in case of standard clustering coefficient yields the following with probability one convergence:  
$$CC^{G\times H} (u,v) \longrightarrow p^2, \quad n \rightarrow \infty.$$

\subsection{$C_2^{G \times H}$ in other random graph models}

In this brief section, we provide a Monte Carlo study for the tensor product of graphs that are sampled from some well-known random graph models. In particular, along with Erd\H{o}s-R\'enyi graphs, we also include the simulation results for random regular graphs \cite{west:01}, Watts-Strogatz small world model \cite{ws:98} and the Bara\'basi-Albert model \cite{ab:02}. Since the current manuscript is mainly on the Erd\H{o}s-R\'enyi model, we do not go into details of the other models, referring to the cited papers for the relevant definitions and theoretical background. 

In order to have a consistent comparison among different models, we fixed the edge density in each experimentation to be $0.5$. Such a choice requires certain selections for the underlying parameters of the random graph models. Our experiments are done in the NetworkX module of Python and the related parameter selections can be found in \cite{eiy:24}. We do not include any further details here in order to avoid repetition. 

Now, the following demonstrates the growth of $C_2^{G \times H}$ when $G$ and $H$ are independently sampled from the stated models:

\begin{figure}[H]
\begin{center}
\includegraphics[scale=0.24]{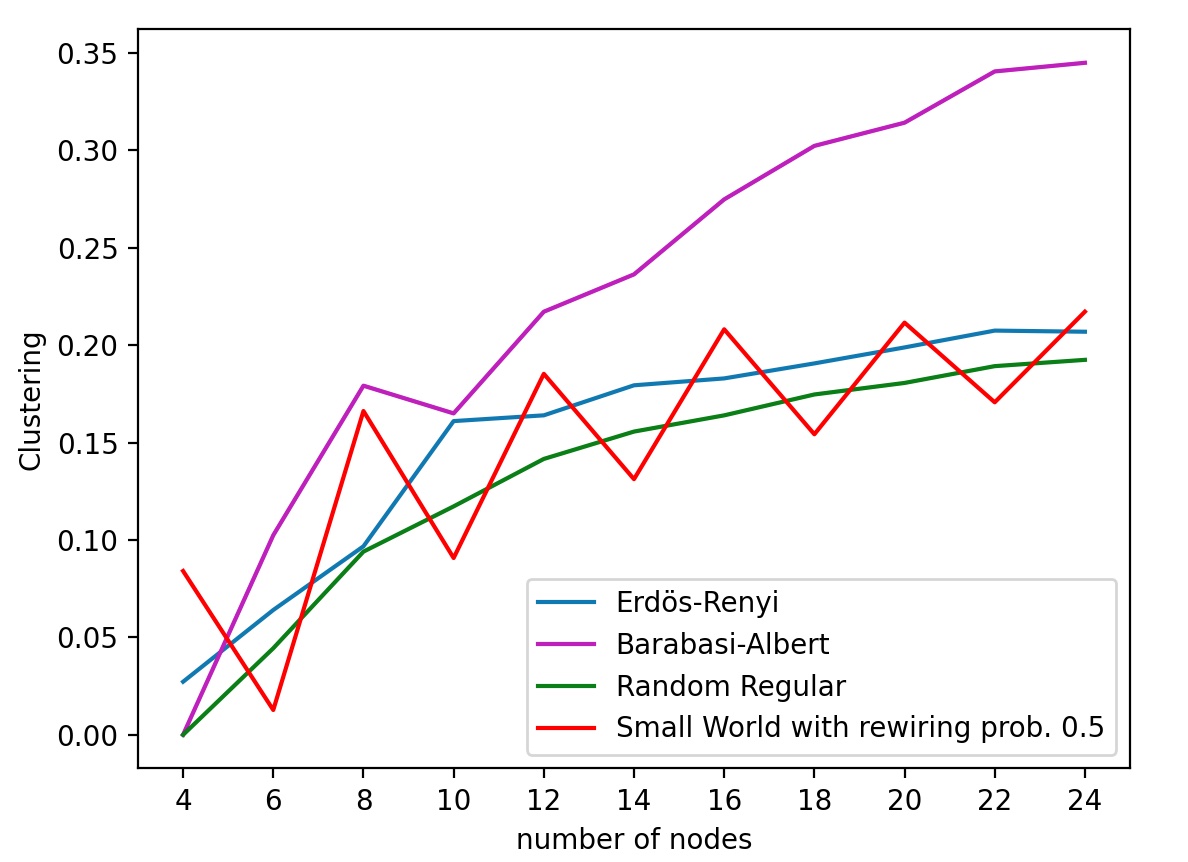}
\end{center}
\end{figure}
\vspace{-0.2in} 

In figure above,  Erd\H{o}s-R\'enyi graph behaves as expected from our theoretical discussion, and the clustering coefficient in this case will approach the value 0.25 as the number of vertices increases. The behaviour in cases of random regular and small world graphs are close to the Erd\H{o}s-R\'enyi case. However, the Bara\'basi-Albert model seems to dominate all others significantly. 
Let us also note that although we have the decomposition given in Proposition \ref{propn:tridecomp},  it is non-trivial to find the expected clustering coefficient exactly for an individual sample from   the model of interest.

\section{Local efficiency}\label{sec:loceff}

We define the  \textit{local efficiency} of a vertex $i$ in some graph $F$  to be $$\eff^F(i) = \frac{1}{\binom{d^K(i)}{2} } \sum_{ \{j,  k \} \subset N(i)} \frac{1}{\dist(j,k)},$$ where $\dist(j,k)$ is the standard graph distance between vertices $j$ and $k$, and $N(i)$ is as before the neighborhood of $i$. (In particular, $N(i)$ does not contain $i$.) Note that for $j, k \in N(i)$, $\dist(j,k)$ is either 1 or 2. 

Our goal in this section is to analyze the local efficiency of a vertex in tensor product of two independent Erd\H{o}s-R\'enyi graphs. For this purpose, we will begin with a general observation relating the local efficiency to the standard clustering coefficient. 

\begin{proposition}\label{propn:eff}
Let $G$ and $H$ be two graphs and $u, v$ be vertices of $G$, $H$, respectively. Then 
$$\eff^{G \times H}((u,v)) = \frac{1}{2} \left(1 + C_2^{G \times H} ((u,v)) \right),$$ where $ C_2^{G \times H} ((u,v))$ is the standard clustering coefficient of $(u,v)$ in $G \times H$.
\end{proposition}

\textbf{Proof.}  
We have, 
\begin{eqnarray*}
 \eff^{G \times H}((u,v)) &=&   \frac{1}{\binom{d^{G \times H} ((u,v))}{2}} \sum_{\{(u_1, v_1), (u_2, v_2) \} \subset N((u,v))} \frac{1}{\dist((u_1, v_1), (u_2, v_2))} \\
 &=&  \frac{1}{\binom{d^{G \times H} ((u,v))}{2}} \sum_{\dist((u_1, v_1), (u_2, v_2)) =1} \frac{1}{\dist((u_1, v_1), (u_2, v_2))} \\ 
 && \hspace{0.2in}  + \hspace{0.1in} \frac{1}{\binom{d^{G \times H} ((u,v))}{2}} \sum_{\dist((u_1, v_1), (u_2, v_2)) =2} \frac{1}{\dist((u_1, v_1), (u_2, v_2))} \\
 &=& \frac{1}{\binom{d^{G \times H} ((u,v))}{2}}  \sum_{\dist((u_1, v_1), (u_2, v_2)) =1} 1 + \frac{1}{\binom{d^{G \times H} ((u,v))}{2}} \sum_{\dist((u_1, v_1), (u_2, v_2)) =2}  \frac{1}{2} \\
 &=& \frac{1}{\binom{d^{G \times H} ((u,v))}{2}}  \sum_{\dist((u_1, v_1), (u_2, v_2)) =1} 1 \\
 && \hspace{0.2in}  + \hspace{0.1in} \frac{1}{2} \frac{1}{\binom{d^{G \times H} ((u,v))}{2}} \left(\binom{d^{G \times H} ((u,v))}{2} - \sum_{\dist((u_1, v_1), (u_2, v_2)) =1} 1  \right) \\
 &=& \frac{1}{2} + \frac{1}{2}\frac{1}{\binom{d^{G \times H} ((u,v))}{2}} \sum_{\dist((u_1, v_1), (u_2, v_2)) =1} 1 =  \frac{1}{2} + \frac{1}{2} \frac{A_2^{G \times H} ((u,v))}{\binom{d^{G \times H} ((u,v))}{2}} \\&=& \frac{1}{2} + \frac{1}{2} C_2^{G \times H} ((u,v)) = \frac{1}{2} \left(1 + C_2^{G \times H} ((u,v)) \right). 
\end{eqnarray*}
\hfill $\square$

Now since by Theorem \ref{thm:asconvofck} we know that $C_2^{G \times H} (u,v) \rightarrow p^2$ with probability one,  Proposition \ref{propn:eff} gives the following result.  

\begin{theorem}
    Let $G$  and $H$ be two independent Erd\H{o}s-R\'enyi  graphs, 
each with parameters $n\in \mathbb{N}$ and $p \in (0,1)$. Let $u \in V(G) $ and $v \in V(H)$. Then, $$\eff^{G \times H}((u,v)) \rightarrow \frac{1}{2} + \frac{1}{2} p^2,$$ with probability one.
\end{theorem}

\section{Isolated vertices}\label{sec:isol}

In this section we are interested in the number of isolated vertices in tensor product of two random graphs. We denote the number of isolated vertices in a given graph $K$ by $I(K)$. The following proposition is elementary and its proof is omitted. 

\begin{proposition}\label{prop:isolatedmoments}
    Let $G$ be an  Erd\H{o}s-R\'enyi graph with parameters $n \in \mathbb{N}$ and $p \in (0,1)$. Then, $$\mathbb{E}[I(G) ] = n (1 - p)^{n-1},$$ 
    and,
    $$\var(I(G)) = n (1 - p)^{n-1} (1 + (np -1) (1 - p)^{n-2}) .$$
\end{proposition}

\begin{theorem}
Let $G$ and $H$ be two independent Erd\H{o}s-R\'enyi graphs with parameters $n \in \mathbb{N}$ and $p \in (0,1)$. We have 
$$\frac{\mathbb{E}[I(G \times H)]}{2 n^2 (1 - p)^{n -1}} \rightarrow 1, \quad n \rightarrow \infty.$$
\end{theorem}

\textbf{Proof.} 
Observe that $(u,v)$ is isolated in $G \times H$ if and only if $u$ is isolated in $G$ or $v$ is isolated in $H$. Then, 
\begin{eqnarray*}
    I(G \times H) &=& \sum_{(u,v) \in V(G \times H)} \mathbf{1}((u,v) \text{ is isolated in } G \times H) \\
    &=&  \sum_{(u,v) \in V(G \times H)} (\mathbf{1}(u \text{ is isolated in } G)  + \mathbf{1}(v \text{ is isolated in } H) \\ && \qquad \qquad \qquad  - \mathbf{1}(u \text{ is isolated in } G, v \text{ is isolated in } H) )\\
    &=&  n I (G) + n I(H) - I (G) I (H).
\end{eqnarray*}
Now, recalling   $\mathbb{E}[I(G)] = n (1 - p )^{n - 1}$, and noting that $G$ and $H$ are independent,
$$\frac{\mathbb{E}[I(G \times H)]}{2 n^2 (1 - p)^{n -1}} = \frac{ 2 n \mathbb{E}[I(G)] - \mathbb{E}[I(G)]   \mathbb{E}[I(H)]}{2 n^2 (1 - p)^{n - 1}} 
   = \frac{2 n^2 (1 - p)^{n - 1} - n^2 (1 - p)^{2 n - 2}}{2n^2 (1 - p)^{n - 1}} \rightarrow 1,  $$
as $n \rightarrow \infty$.  \hfill $\square$

The theorem in particular says that  $\mathbb{E}[I(G \times H)] \approx 2 n^2 (1 - p)^{n - 1} $ for large $n$. 

\begin{remark}
    The criteria for being an isolated vertex, i.e. being of degree $0$, given in the proof of previous theorem can be generalized to  vertices in $G \times H$ with degree $m\in \mathbb{N}$. Namely, if we denote the number of degree $m$ vertices in a graph $K$ by $D_m^K$, then, 
    \begin{eqnarray*}
    D_m^{G \times H} = \sum_{(u,v) \in V(G \times H)} \mathbf{1}(d^G(u) = d^H(v) = m) &+& \sum_{(u,v) \in V(G \times H)} \mathbf{1}(d^G(u) = m, d^H(v) > m) \\ && + \sum_{(u,v) \in V(G \times H)} \mathbf{1}(d^G(u) > m, d^H(v) = m).    
    \end{eqnarray*} 
    We will not go into further details for this more general case here. 
\end{remark}

Focusing on the number of isolated vertices, when $p = c / n$, for some $c > 1$, we will  further prove  a weak law of large numbers. This result will require the following proposition whose proof is again elementary and is therefore omitted. 

\begin{proposition}\label{propn:eltvar}
    Let $X, Y $ be independent and identically distributed random variables with finite variance. Then 
    $$\var (X Y ) = \var (X) (\var (X)  + 2 (\mathbb{E}[X])^2).$$
\end{proposition}

\begin{theorem}
  Let $G$ and $H$ be two independent Erd\H{o}s-R\'enyi graphs with parameters  $n \in \mathbb{N}$ and $p \in (0,1)$. We assume that $p = c /n$ for some positive constant $c  > 1$.  Then we have $$\frac{I(G \times H)}{2n^2 (1 - p)^{n-1} - n^2 (1-p)^{2 n -2}} \rightarrow_{\mathbb{P}} 1,$$ as $n \rightarrow \infty$, where the convergence is in  probability. 
\end{theorem} 

\textbf{Proof.} Note here that $\mathbb{E}[I(G  \times H)] = 2n^2 (1 - p)^{n-1} - n^2 (1-p)^{2 n -2}$. Let $\epsilon > 0$. Recalling  $I(G \times H) = n I (G) + n I(H) - I (G) I (H)$, we have, 
\begin{eqnarray*}
    \mathbb{P} \left( \left| \frac{I(G \times H) - \mathbb{E}[I(G  \times H)]}{\mathbb{E}[I(G  \times H)] } \right|  > \epsilon \right) &\leq&  \mathbb{P} \left( \left| \frac{n (I(G)  + I(H)) - \mathbb{E}[n (I(G)  + I(H))}{\mathbb{E}[I(G  \times H)] } \right|  > \frac{\epsilon}{2} \right) \\
    && +  \mathbb{P} \left( \left| \frac{I(G) I(H)  - \mathbb{E}[I(G) I(H) ]}{\mathbb{E}[I(G  \times H)] } \right|  > \frac{\epsilon}{2} \right) \\
    &\leq& \frac{\var(n(I(G) + I(H)))}{(\mathbb{E}[I(G  \times H)])^2 } \frac{4}{\epsilon^2} +   \frac{\var(I (G) I (H))}{(\mathbb{E}[I(G  \times H)])^2 } \frac{4}{\epsilon^2}  
\end{eqnarray*}
Now,  looking at the first term on the right-most when $p = c/n$, using Proposition \ref{prop:isolatedmoments}, we have $$\frac{\var(n(I(G) + I(H)))}{(\mathbb{E}[I(G  \times H)])^2 } = \frac{2n^2  n (1 - p)^{n-1} (1 + (np -1) (1 - p)^{n-2})}{(2n^2 (1 - p)^{n-1} - n^2 (1-p)^{2 n -2})^2} = \mathcal{O} \left(\frac{1}{n} \right).$$ 
For the second term on the right-most,  we use Proposition \ref{propn:eltvar} to get, 
$$ \frac{\var(I (G) I (H))}{(\mathbb{E}[I(G  \times H)])^2 }   = \frac{\var (I(G)) (\var (I(G))  + 2 (\mathbb{E}[I(G)])^2)}{(\mathbb{E}[I(G  \times H)])^2 }.$$ 
Now, again using Proposition \ref{prop:isolatedmoments} and doing  some straightforward  computations give $$ \frac{\var(I (G) I (H))}{(\mathbb{E}[I(G  \times H)])^2 }  = \mathcal{O} \left(\frac{1}{n} \right).$$
Result follows. \hfill $\square$

\section{Conclusion}\label{sec:conc}

In this paper, we studied the tensor product of independent Erd\H{o}s-R\'enyi graphs. In particular, when $p=1/2$ we showed that largest clique in the product has size  approximately $2 \log n$ as in the case of an individual Erd\H{o}s-R\'enyi graph. Also, using techniques similar to the analysis of the largest clique, we obtained a strong law of large numbers for an extended clustering coefficient. Other contributions included observations on the local efficiency and the number of isolated vertices in the tensor product.  

However, as already noted in the Introduction, our initial motivation was on studying the  modular product of two independent Erd\H{o}s-R\'enyi graphs, following the fact that  the size of the maximum common subgraph of two graphs $G$ and $H$ turns out to be the size of the largest clique in modular product of $G$ and $H$. Although, certain parts of our calculations on cliques carry over to the case of modular products, we are currently not able to fully understand the mean of the largest clique in modular product case. We hope to pursue in this in an upcoming work. That said, let us also include some simulation results for the size of the  largest clique in the modular product for small values of $n$: 

\begin{figure}[H]
\begin{center}
\includegraphics[scale=0.3]{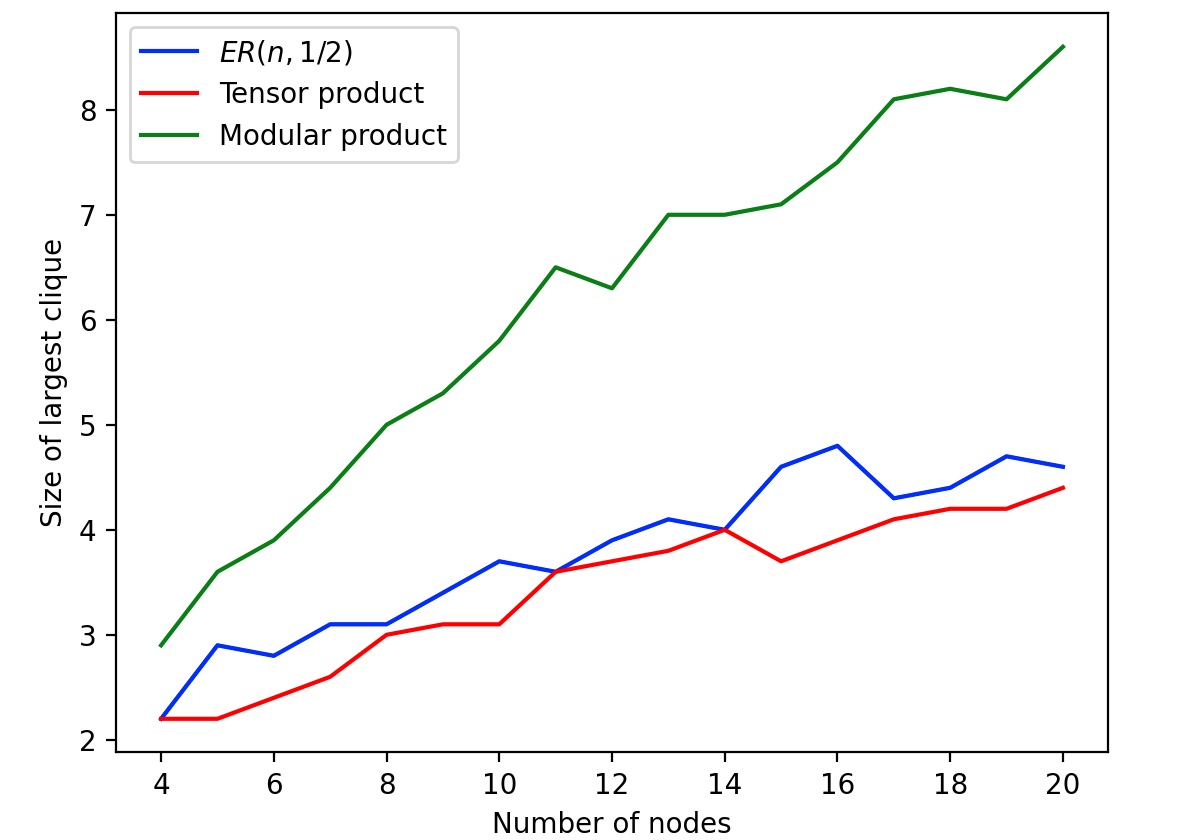}
\end{center}
\end{figure}

\vspace{-0.15in}

The figure demonstrates average maximal clique size in three settings: (1) in a single Erd\H{o}s-R\'enyi graph, (2) tensor product of two independent Erd\H{o}s-R\'enyi graphs, and (3) modular product of two independent Erd\H{o}s-R\'enyi graphs. In this Monte Carlo study, for each of the three cases,  25 samples are used  in order to compute the average maximal clique size. 

The behaviours for the former two agree as we already know from the theory part. On the other hand, the size of the maximal clique in the corresponding modular product seems to behave twice as the other two in tensor product setting. This is intuitively clear once one  observes that the probability of presence for any edge is twice as much as the corresponding probability in tensor product case. However, the computations we did for the tensor product case do not work in here, and as already noted, we plan to analyze the modular product case in more detail in a subsequent work.  

As a last note, and as a second possible research direction, it would be interesting to understand graph products of random graphs for models other than the Erd\H{o}s-R\'enyi model. Certain generalizations to inhomogeneous random graphs seem to be at reach, but one may also consider other models from the literature such as the small world model of Watts and Strogatz and the preferential attachment model of Albert and Barab\'asi.

\section*{Appendix}

\subsection*{Proof of Proposition \ref{propn:basicupbound}}

 For $A \subset \{1,2,\ldots,n\}$ with $|A|= k$, let $I_A$ be the indicator that the clique corresponding to $A$ is present in $G$. Then, 
$$\var(X_k) = \var \left(\sum_{A} I_A \right) = \sum_{A} \var (A) + \sum_{A \neq B} \cov(I_A, I_B).$$ (Sets in indices are understood to be of size $k$.)
Observe that the covariance have non-zero effect only when  $A, B$ contain at least two common vertices. We denote the set of such $(A,B)$ pairs by  $\D$. Using the obvious estimates $\var(I_A) \leq \mathbb{E}[I_A^2] = \mathbb{E}[I_A] $ and $\cov(I_A, I_B) \leq \mathbb{E}[I_A I_B]$, we then have $$\var(X_k) \leq \sum_{A} \E[I_A] + \sum_{(A,B) \in \D } \E[I_A I_B] = \E[X_k] +  \sum_{(A,B) \in \D } \E[I_A I_B].$$
 Focusing on $\sum_{(A,B) \in \D } \E[I_A I_B]$, we have 
 \begin{eqnarray*}
    \sum_{(A,B) \in \D } \E[I_A I_B] &=& \sum_{(A,B) \in \D }\mathbb{P}(I_A = 1, I_B = 1) =  \sum_{(A,B) \in \D }\mathbb{P}(I_A = 1)  \mathbb{P} (I_B = 1 \mid  I_A = 1) \\
    &=& \sum_A \sum_{B: (A,B) \in \D } \mathbb{P}(I_A = 1)  (I_B = 1 \mid  I_A = 1) \\ &=&  \sum_A  \mathbb{P}(I_A = 1) \sum_{B: (A,B) \in \D }  \mathbb{P}  (I_B = 1 \mid  I_A = 1).
 \end{eqnarray*}
Now note that $\sum_{B: (A,B) \in \D } \mathbb{P}(I_A = 1)  (I_B = 1 \mid  I_A = 1)$ is independent of $A$ due to symmetry. So let us choose  and fix some $A_* \subset \{1,2,\ldots,n\}$ of size $k$ and write $\sum_{B: (A_*,B) \in \D }  (I_B = 1 \mid  I_{A_*} = 1)$ for this common probability. Then, $$ \sum_{(A,B) \in \D } \E[I_A I_B] =  \sum_A  \mathbb{P}(I_A = 1) \sum_{B: (A,B) \in \D }  \mathbb{P}  (I_B = 1 \mid  I_A = 1) = \mathbb{E}[X_k]\sum_{B: (A_*,B) \in \D }  \mathbb{P}  (I_B = 1 \mid  I_{A_*} = 1). $$
Now, we observe that $$\sum_{B: (A_*,B) \in \D }  \mathbb{P}  (I_B = 1 \mid  I_{A_*} = 1) = \sum_{i=2}^{k-1} \sum_{B: (A_*,B) \in \D, |A_* \cap B| = i }  \mathbb{P}  (I_B = 1 \mid  I_{A_*} = 1) = \sum_{i=2}^{k-1} \binom{k}{i} \binom{n-k}{k-i}p^{(\binom{k}{2} - \binom{i}{2} )}. $$
Result follows by combining the observations we have made. \hfill $\square$

\subsection*{Proof of Corollary \ref{cor:varbound}.}

(i)
We already know that 
\begin{eqnarray*}
\var(X_k) &\leq&  \binom{n}{k}  p^{\binom{k}{2}} + \binom{n}{k}  p^{\binom{k}{2}} \sum_{i=2}^{k-1} \binom{k}{i} \binom{n  - k}{k - i} p^{\binom{k}{2} - \binom{i}{2}} \\
&=& \binom{n}{k}  p^{\binom{k}{2}} + \binom{n}{k}  p^{\binom{k}{2}} p^{\binom{k}{2}} \sum_{i=2}^{k-1} \binom{k}{i} \binom{n  - k}{k - i} p^{ - \binom{i}{2}}     \\
&\leq&  \binom{n}{k}  p^{\binom{k}{2}} + \binom{n}{k}  p^{k(k-1)} \sum_{i=2}^{k-1} \binom{k}{i} \binom{n  - k}{k - i} p^{ - \binom{i}{2}}.   
\end{eqnarray*}
Now we are willing to understand $\sum_{i=2}^{k-1} \binom{k}{i} \binom{n  - k}{k - i} p^{ - \binom{i}{2}}$.

Let $a_i = \binom{k}{i} \binom{n  - k}{k - i} p^{- \binom{i}{2}}$. We claim that $a_i$ is decreasing for $i = 2,3,\ldots,k-1$ when $n$ is sufficiently large. We have,
\begin{equation}\label{eqn:ratio1}
\frac{a_{i+1}}{a_i} = \frac{k - i}{ n - 2k + i + 1} \frac{k - i }{ i + 1} p^{-i} \leq  \frac{1}{n - 2k + i + 1} \frac{k^2}{i+1} p^{-i}.    
\end{equation}
When $i =2$, $$\frac{a_3}{a_2} \leq \frac{k^2 p^{-2}}{3(n -2k +3)},$$ and the right-hand side is clearly less than 1 when $n $ is large, and $k$ is as in the statement of the corollary.

Now, for $i \in \{3,4,\ldots,k-2\}$,  observe that $$ \frac{1}{n - 2k + i + 1}  - \frac{1}{n} = \frac{2 k - i -1}{(n - 2k + i + 1)n } \leq \frac{2k }{(n - 2k + i + 1) n} \leq \frac{4k}{n^2},$$ where in the last step we used $n - 2 k+ i + 1 \geq n/2$ for large enough $n$, and for $k$ given in statement of the corollary. This says $$\frac{1}{n - 2k + i + 1} \leq \frac{1}{n} + \frac{4k}{n^2}.$$ Using this in \eqref{eqn:ratio1} for $i \in \{3,4,\ldots,k-2\}$, we obtain, 
$$\frac{a_{i+1}}{a_i} \leq \left( \frac{1}{n} + \frac{4k }{n^2}\right) \frac{k^2}{i+1}p^{-i}.$$
Then, using $i! \geq (i/e)^i >   (i/e)^{i-2}$ for $i \geq 3$, $e^x \geq 1 + x$ for all $x$, and  $e^{4k^2 /n }\leq 2$ for large $n$, we obtain, 
\begin{eqnarray*}
    \frac{a_{i}}{a_2} &=& \prod_{j=2}^{i-1} \frac{a_{j+1}}{a_j} \leq \prod_{j=2}^{i-1} \left( \frac{1}{n} + \frac{4k }{n^2}\right) \frac{k^2}{j+1}p^{-j} = \frac{1}{i!} \left( \frac{k^2}{n}\right)^{i-2} p^{- \frac{(i - 2)(i+1)}{2}} \prod_{j=2}^{i-1} \left(1 + \frac{4k}{n} \right) \\
    &\leq&    \left(\frac{e}{i} \right)^{i-2}  \left( \frac{k^2}{n}\right)^{i-2} \left(p^{- (i+1)/2} \right)^{i -2} \left(e^{\frac{4k}{n}}\right)^k \leq 2 \left( \frac{ek^2}{in } p^{-(i+1)/2}\right)^{i - 2}. 
\end{eqnarray*}
We claim that $\frac{ek^2}{in } p^{-(i+1)/2} < \frac{1}{2}$ for $i \in \{3,4,\ldots,k-2\}$ for large enough $n$. We have, 
\begin{eqnarray*}
    \frac{ek^2}{in } p^{-(i+1)/2} &\leq& \left( \frac{e\sqrt{p^{-1}}}{3} \right) \frac{k^2}{n} p^{-k/2} = \left( \frac{e\sqrt{p^{-1}}}{3} \right) \frac{k^2}{n}  p^{-\frac{2 \log_{1/p}n - M \log_{1/p} \log n}{2}} \\
    &=& \left( \frac{e\sqrt{p^{-1}}}{3} \right) \frac{k^2}{n} p^{\log_p n  - \frac{M}{2} \log_p \log n} = \left( \frac{e\sqrt{p^{-1}}}{3} \right) \frac{k^2}{n}  n (\log n)^{- M/2} \\
    &=&  \left( \frac{e\sqrt{p^{-1}}}{3} \right)  k^2  (\log n)^{- M/2}.
\end{eqnarray*}
Now, for $M > 4$, $$\frac{k^2}{(\log n)^{M/2}} \leq \frac{2 (\log_{1/p} n)^2}{(\log n)^{M/2}} \rightarrow 0,$$ as $n \rightarrow \infty$. In particular, for $i \in \{3,4,\ldots,k-2\}$,  $\frac{ek^2}{in } p^{-(i+1)/2}  < \frac{1}{2}$ for large enough $n$. Thus the sequence $a_i$ is decreasing for given selection of $k$, when $n$ is large enough.  Hence, we may use $a_2$ as the corresponding maximum value. 
But then we obtain 
\begin{eqnarray*}
\var(X_k)  
&\leq&  \binom{n}{k}  p^{\binom{k}{2}} + \binom{n}{k}  p^{k(k-1))} \sum_{i=2}^{k-1} \binom{k}{i} \binom{n  - k}{k - i} p^{ - \binom{i}{2}} \\
&\leq& \binom{n}{k}  p^{\binom{k}{2}} + \binom{n}{k}  p^{k(k-1))} k \binom{k}{2} \binom{n-k}{k-2} p^{-1} \\
&\leq& \binom{n}{k}  p^{\binom{k}{2}} + \binom{n}{k}  p^{k(k-1))} \frac{k^3}{2}  \binom{n-k}{k-2}  p^{-1},
\end{eqnarray*}
as required. 

(ii) The proof follows similarly for fixed $k$ case. In particular, the key argument above was showing that the sequence $a_i$ is decreasing over $2,3,\ldots,k-2$. Here showing that $a_3 < a_2 $ for large enough $n$ is  exactly the same as above. For $i \in \{3,4,\ldots,k-2\}$, we still have,
$$\frac{a_i}{a_2} \leq 2 \left( \frac{ek^2}{in } p^{-(i+1)/2}\right)^{i - 2},$$ but this time the remaining part is easier since $$\frac{ek^2}{in } p^{-(i+1)/2} \leq \frac{ek^2}{3n } p^{-(k+1)/2},$$ and the right-hand term clearly converges to 0 since $k$ is fixed. The remaining parts are treated as in (i).
\hfill $\square$

\vspace{0.1in} 

\textbf{Acknowledgements.} First author has been supported by the Scientific and Technological Research
Council of Turkey TUBITAK Grant No 122M452, and  BAGEP Award of the Science Academy, Turkey.

\end{document}